\hfill {\bf IC/97/92} 
\baselineskip=16pt
\font\twelvebf=cmbx12

\vskip 15mm
\def\Z{{\bf Z}}
\def\C{{\bf C}}
\def\Ch{{\C[[h]]}}
\def\N{{\bf N}}
\def\Q{{\bf Q}}
\def\e{\eqno}
\def\l{\ldots}
\def\n{\noindent}
\def\UA{U_h(A_\infty)}
\def\Ua{U_h(a_\infty)}
\def\gl{U_h(gl_\infty)}
\def\q{\quad}
\def\V{V(\{M\};\xi_0,\xi_1)}
\def\G{\Gamma(\{M\})}
\def\a{\alpha}
\def\b{\beta}
\def\g{\gamma}
\def\r{\rho}
\def\Vf{V(\{M\}_{m,n};M_m,M_n)}
\def\t{\theta}
\def\Vh{{\hat V}(\{M\}_{m,n};M_m,M_n)}

\n
{\twelvebf Highest weight irreducible representations of the
quantum algebra
U$_h$(A$_\infty$)}

\vskip 32pt
\noindent
T.D. Palev\footnote*{Permanent address: Institute for Nuclear
Research and Nuclear Energy, 1784 Sofia, Bulgaria; E-mails:
tpalev@inrne.acad.bg, stoilova@inrne.acad.bg}

\noindent
International Centre for Theoretical Physics, 34100 Trieste, Italy
\vskip 12pt

\vskip 12pt

\noindent
N.I. Stoilova*

\noindent
Department of Mathematics, University of Queensland, Brisbane Qld
4072,
Australia

\vskip 4cm
\noindent
Running title: Representations of the quantum algebra $U_h(A_\infty)
$

\vskip 48pt
\noindent
{\bf Abstract.} A class of highest weight irreducible
representations of the algebra $\UA$, the quantum analogue of the
completion and central extension $A_\infty$ of the Lie algebra
$gl_\infty$, is constructed. It is considerably larger than the
known so far representations. Within each module a basis is
introduced and the transformation relations of the basis under
the action of the Chevalley generators are explicitly written.
The verification of the quantum algebra relations to be satisfied
is shown to reduce to a set of nontrivial $q$-number identities.
All our representations are restricted in the terminology of S.
Levendorskii and Y. Soibelman (Commun. Math.  Phys. {\bf 140},
399-414 (1991)).

\vskip 4cm

\vfill\eject
\leftskip 0pt
\vskip 12pt

\n
{\bf 1. INTRODUCTION}

\bigskip\n
We construct a class of highest weight irreducible
representations (irreps) of the algebra $\UA \;[8],$ the quantum
analogue of the completion and central extension $A_\infty$ of
the Lie algebra $gl_\infty$ $[1],[6].$ Our interest in the subject
stems from the observation that certain representations of
$gl(n)$ (including $n=\infty$) $[9]$ and of $A_\infty$ (see {\it
Example 2} in Ref. [10] and the references therein) are related to a
new quantum statistics, the $A$-statistics. The latter, as it is
clear now, is a particular case of the Haldane exclusion
statistics $[5],$ a subject of considerable interest in condensed
matter physics (see Sect. 4 in Ref. [13] for more discussions on the
subject). It turns out that some of the representations of the
deformed algebra $\UA$ satisfy also the requirements of the
Haldane statistics. More precisely, they lead to new solutions
for the microscopic statistics of Karabali and Nair $[7]$
directly in the case of infinitely many degrees of freedom. These
results will be published elsewhere. We mention them here only in
order to justify our personal motivation for the work we are
going to present. One may expect certainly that similar as for
$A_\infty$ $[2], [4]$ the representations of $\UA$ may prove useful
also in other branches of physics and mathematics.

The quantum analogues of $gl_\infty$ and $A_{\infty}$ in the
sense of Drinfeld $[3],$ namely $U_h(gl_\infty)$ and
$U_h(A_\infty)$, were worked out by Levendorskii and Soibelman $[8].$
These authors have constructed  a class of highest
weight irreducible representations, writing down explicit
expressions for the transformations of the basis under the action
of the algebra generators.

The $U_h(A_\infty)$-modules, which we study, are labeled by all
possible sequences (see the end of the introduction for
the notation)
$
\{M\} \equiv \{M_i \}_{i\in \Z}\in \C[h]^\infty
$,
subject to the conditions:

(a) There exists $m\le n \in \Z$, such that $M_m=M_{m-k}$ and
   $M_n=M_{n+k}$ for all $k\in \N$;

(b) $M_i-M_j\in \Z_+ $ for all $i<j \in \Z.$

\smallskip\n
Representations, corresponding to two different sequences,
$\{M^1\} \ne \{M^2\}$, are inequivalent.  The
$U_h(A_\infty)$-modules of Levendorskii and Soibelman $[8]$
consist of all those sequences $\{M^{(s)}\}$, for which $s=m=n\in
\Z$ and $M_i^{(s)}=1$, if $i< s$ and $M_i^{(s)}=0$ for $i\ge s$.

In Refs. [11] and [12] a class of highest weight irreps, called finite
signature representations, of the Lie algebra $A_\infty$ was
constructed. The corresponding modules are 
labeled by the set
of all complex sequences
$
\{M\} \equiv \{M_i \}_{i\in \Z}\in \C^\infty
$,
which satisfy the conditions (a) and (b). The name "finite
signature" comes to indicate that, due to (a), each signature
$\{M\}$ is characterized by a finite number of different
coordinates or, more precisely, by no more then $n-m+1$ different
complex numbers.

From our results it follows 
that each $A_\infty-$module with a signature
$\{M\}$ can be deformed to an $\UA-$module with the same
signature. The class of the finite signature representations of
$\UA$ is however larger, which is due to the fact that the
coordinates of the $\UA$ signatures take values in $\C[h]$.  All
representations we obtain are restricted in the sense of Ref.
[8], Definition 4.1.

In Section 2 we construct a class of highest weight irreps of
the subalgebra $\Ua$ of $\UA$. Some of these representations,
namely the finite signature representations ({\it Definition 2}),
are extended to representations of $\UA$  in Section 3.

Throughout the paper we use the notation
(most of them standard):

$\N$ - all positive integers;

$\Z_+$ - all non-negative integers;

$\Z$ - all integers;

$\Q$ - all rational numbers;

$\C$ - all complex numbers;

$\C[h]$ - the ring of all polynomials in $h$ over $\C$;

$\C[[h]]$ - the ring of all formal power series in $h$ over $\C$;

$X^n = \{(a_1,a_2,\l,a_n) | \; a_i\in X  \}$, (including $n=\infty$);

$[a,b]=\{x|\; a\le x \le b, x\in \Z \},\q (a,b)=\{x|\; a < x < b, x\in \Z \}$;

$\Gamma(\{M\})$ - the $C$-basis of a module with a signature
$\{M\}$;

$
\theta (i) = \cases
{ 1, & for $i\ge 0$ \cr
  0, &  for $i<0$; \cr}
$

$q=e^{h/2} \in \Ch$;

$
[x]={q^x-q^{-x}\over {q-q^{-1}}}\in \C[[h]];
$

$[x,y]_q=xy-qyx.$

\vskip 24pt

\n
{\bf 2. REPRESENTATIONS OF THE ALGEBRA U$_h$($a_\infty$)}

\bigskip\n
The algebra $\Ua$ was defined in Ref. [8]. The authors denote it
as $U_h(g'(A_\infty))_f$. It is a Hopf algebra, which is a
topologically free module over $\C[[h]]$ (complete in $h$-adic
topology), with generators $\{e_i, f_i, h_i, c  \}_{i\in \Z}$,
and
\smallskip
\noindent
1. Cartan relations:

$$
\eqalign{
& [c, a]=0, \quad a\in \{h_i, e_i, f_i \}_{i\in \Z} \cr
& [h_i, h_j]=0, \cr
& [h_i, e_j]=(\delta _{ij}-\delta _{i,j+1})e_j,  \cr
& [h_i, f_j]=-(\delta _{ij}-\delta _{i,j+1})f_j,  \cr
& [e_i, f_j]=\delta_{ij}{{q^{h_i-h_{i+1}+(\theta (-i)-\theta (-i-1))c}
-q^{-h_i+h_{i+1}-(\theta (-i)-\theta (-i-1))c}}
  \over{q-q^{-1}}}.\cr
}\e(1)
$$
\noindent
2. $e$-Serre relations:

$$
\eqalign{
& e_ie_j=e_je_i, \quad {\rm if} \;\; |i-j|\neq 1,\cr
& e_i^2e_{i+1}-(q+q^{-1})e_ie_{i+1}e_i+e_{i+1}e_i^2=0, \cr
& e_{i+1}^2e_i-(q+q^{-1})e_{i+1}e_ie_{i+1}+e_ie_{i+1}^2=0. \cr
}\e(2)
$$

\n
3. $f$-Serre relations:

$$
\eqalign{
& f_if_j=f_jf_i, \quad {\rm if} \;\; |i-j|\neq 1,\cr
& f_i^2f_{i+1}-(q+q^{-1})f_if_{i+1}f_i+f_{i+1}f_i^2=0, \cr
& f_{i+1}^2f_i-(q+q^{-1})f_{i+1}f_if_{i+1}+f_if_{i+1}^2=0. \cr
}\e(3)
$$

\bigskip
We do not write the other Hopf algebra maps ($\Delta,\;
\varepsilon,\; S$) $[8],$ since we will not use them. They
are certainly also a part of the definition
of $\Ua$.

Replacing throughout in the above relations
$\{e_i, f_i, h_i, c  \}_{i\in \Z}$ with
$\{E_i, F_i,H_i,0\}_{i\in \Z}$, one obtains the definition of
$U_h(gl_\infty)$.

\bigskip
In terms of an equivalent set of generating elements
$\{ {\hat e}_{i},{\hat f}_{i}, h_i,c \}_{i\in \Z}$, with
$$
{\hat e}_{i}=e_i q^{(h_{i+1}-h_i)/2},\q
{\hat f}_{i}=f_i q^{(h_{i}-h_{i+1})/2},\e(4)
$$
one writes the quantum analogues of the Weyl generators
$\{e_{ij}\}_{(i,j)\in \Z^2}$:
$$
\eqalign{
& e_{ii}=h_i,\q e_{i,i+1}={\hat e}_{i},\q e_{i+1,i}={\hat f}_{i},  \cr
& e_{ij}= [{\hat e}_{i},[{\hat e}_{i+1},[\l,
[{\hat e}_{j-2},{\hat e}_{j-1}]_q \l ]_q]_q]_q, \q i+1<j,  \cr
& e_{ji}= [{\hat f}_{i},[{\hat f}_{i+1},[\l,
[{\hat f}_{j-2},{\hat f}_{j-1}]_q \l ]_q]_q]_q, \q i+1<j.  \cr
}\e(5)
$$
The "commutation relations" between these generators follow from
(1)-(3) and are given in Ref. [8]. The relevance of the
generators $\{e_{ij}\}_{(i,j)\in \Z^2}$ stems from the
observation that the set of ordered monomials (see Ref. [8] for
the ordering)
$$
c^l\prod_{(i,j)\in \Z^2}e_{ij}^{n_{ij}} \e(6)
$$
with finitely many non-zero exponents $n_{ij}\in \Z_+,\;l\in \Z_+$
forms a (topological) basis in $\Ua$.

Set $H=\oplus _i \C h_i,$ define linear functionals
$\varepsilon _i:H\rightarrow \C $ by $\varepsilon _i(h_j)=\delta _{ij}$
and set $\alpha _i=\varepsilon _i-\varepsilon _{i+1},\;\;
Q_+'=\oplus _i\Z_+\alpha _i.$ Denote by $U_h(n_+)$ (respectively
$U_h(n_-)$ ) the unital subalgebra in $U_h(a_\infty )$ generated
by $\{ e_i\}_{i\in \Z}$  (respectively $\{ f_i\}_{i\in \Z}$ ).
Then

$$
U_h(n_\pm )=\oplus _{\alpha \in Q_+'}U_h(n_\pm )_{\pm \alpha},
\e (7)
$$
where

$$
U_h(n_\pm )_{\pm \alpha}=\{ x\in U_h(n_\pm ) |\; [h',x]=\pm \alpha (h')x,\;
\forall h'\in H \} \e (8)
$$
for $\alpha \neq 0,$ and $U_h(n_\pm )_0=\Ch.$  Any element
$u\in U_h(a_\infty )$ can be represented as

$$
u=\sum_{k=0}^\infty h^k \sum_{l=0}^{ l(k) } c^l
\sum_{\alpha , \beta \in Q_+'} \; \sum _{\gamma (k,l)\in \Z_+^\infty }
\sum_{t=1}^{t(\alpha,\beta)}
{\cal F}_{\alpha ,k,l,t}\prod_{i\in \Z} h_{i}^{\gamma (k,l)_i}
{\cal E} _{\beta ,k,l,t},   \;\;  finite \; sums \; over
\; \alpha,\; \beta, \; \gamma,\e(9)
$$
where ${\cal F} _{\alpha ,k,l,t}\in U_h(n_- )_{-\alpha},\;
{\cal E}_{\beta ,k,l,t}\in U_h(n_+ )_{+\beta } $
and finitely many exponents $\gamma (k,l)_i$ are different from
zero. The words "$finite \; sums \; over
\; \alpha,\; \beta, \; \gamma$"  have been added in order to
indicate that for a fixed $k$ only finitely many summands in (9) are
different from zero.

The set ${\hat U}_h(a_\infty)$, consisting of all $\Ch-$polynomials
of the Chevalley generators $\{e_i, f_i, h_i, c  \}_{i\in \Z}$,
is dense in $\Ua$ with a basis (6). In particular
${\cal F} _{\alpha ,k,l,t}$, ${\cal E}_{\beta ,k,l,t}$ and
$\prod_{i\in \Z} h_{i}^{\gamma (k,l)_i}$ are in
${\hat U}_h(a_\infty)$. Then according to (9) any element
$u\in \Ua$ is of the form
$$
u=\sum_{k=0}^\infty u_k h^k,\q \q
u_k=\sum_{l=0}^{ l(k) } c^l
\sum_{\alpha , \beta \in Q_+'}\;\sum _{\gamma (k,l)\in \Z_+^\infty }
\sum_{t=1}^{t(\alpha,\beta)}
{\cal F}_{\alpha ,k,l,t}\prod_{i\in \Z} h_{i}^{\gamma (k,l)_i}
{\cal E} _{\beta ,k,l,t}\in {\hat U}_h(a_\infty). \e(10)
$$

We pass to construct a class of highest weight irreps of
$\Ua$.  To this end we define the $\Ua$-module
$V(\{M\};\xi_0,\xi_1)$, 
labeled by $\xi_0,\; \xi_1 \in \C [h]$ and by a
sequence $\{M\}\equiv \{M_i\}_{i\in \Z }\in \C [h]^\infty$ such that
$$
M_i-M_j\in \Z_+,\;\; \forall i<j\in \Z,
$$
with a basis $\G$. The latter,
called a central basis ($C$-bases), is independent on
$\xi_0,\xi_1$.  It is formally the same as the one introduced in
Refs. [11] and [12] for a description of representations of
$a_\infty$. $\G$ consists of all $C-$ patterns

$$
(M)\equiv \left[\matrix
{ .., & M_{1-\theta -k},& \ldots, &M_{-1},&M_0, &M_1, &\ldots,
&M_{k-1},...\cr
.., & \ldots & \ldots & \ldots & \ldots &\ldots  &\ldots &\ldots \cr
&M_{1-\theta -k,2k+\theta -1}, & \ldots, &M_{-1,2k+\theta -1},
&M_{0,2k+\theta -1}, & M_{1,2k+\theta -1}, & \ldots,
& M_{k-1,2k+\theta -1} \cr
&\ldots &\ldots & \ldots & \ldots  &\ldots &\ldots \cr
& & & M_{-1,3}, & M_{03}, &M_{13} \cr
& & & M_{-1,2}, &M_{02} \cr
& & & & M_{01} \cr
}\right], \eqno(11)
$$
where $k\in \N,\; \theta=0,1$.
Each such pattern is an ordered collection of formal polynomials
in $h$
$$
M_{i, 2k+\theta -1}\in \C[h], \quad \forall k \in \N, \quad \theta =0,1,
\quad
\;\; i \in [-\theta -k+1, k-1], \eqno(12)
$$
which satisfy the conditions:
\smallskip
\noindent
($i$) there exists a positive, depending on $(M)$, integer
$N_{(M)} > 1$, such that
$$
M_{i,2k+\t-1}=M_i, \;\; \forall \; 2k+\t-1 \geq N_{(M)}, \quad
 \theta =0,1, \;\; i \in [1-\theta -k, k-1]; \eqno(13a)
$$
(ii) for each $k \in \N,  \; \theta =0,1$ and $ i \in [1-\theta -k, k-1]$
$$
M_{i+\theta -1,2k+\theta }-M_{i, 2k+\theta -1} \in \Z_+ ,
\;
M_{i,2k+\theta -1}-M_{i+\theta , 2k+\theta} \in \Z_+. \eqno(13b)
$$

Denote by ${\hat V}(\{M\};\xi_0,\xi_1)$ the free $\C[[h]]$-module
with generators $\G$ and let $\V$ be its completion in the
$h$-adic topology. $\V$ is a topologically free $\C[[h]]$-module
with a (topological) basis $\G$ and ${\hat V}(\{M\};\xi_0,\xi_1)$
is dense in it (in the $h$-adic topology).
$\V$ consists of all formal power series in $h$ with
coefficients in ${\hat V}(\{M\};\xi_0,\xi_1)$:
$$
v=\sum_{i=0}^\infty v_ih^i,
\quad v_0, v_1, v_2, \ldots \in {\hat V}(\{M\};\xi_0,\xi_1). \e(14)
$$
If $a$ is a $\Ch-$linear map in ${\hat V}(\{M\};\xi_0,\xi_1)$,
$a\in End\; {\hat V}(\{M\};\xi_0,\xi_1)$, we extend it
to a continuous linear map on $\V$ setting
$$
av=\sum_{i=0}^\infty (av_i)h^i.\e(15)
$$
Therefore the transformation of $\V$ under the action of
$a$ is completely defined, if  $a$ is defined on
$\G$.

We proceed to turn $\V$ into a  $\Ua$ module. Denote by $(M)_{\pm
\{ j,p\} }$ and $(M)_{\pm \{ j,p\} }^{\pm \{ l,q\} }$ the
patterns  obtained from the $C$-pattern $(M)$ in (11) after the
replacements
$$
M_{jp}\rightarrow M_{jp}\pm 1 \q {\rm and} \q
M_{lq}\rightarrow M_{lq}\pm 1, \quad
M_{jp}\rightarrow M_{jp}\pm 1,
$$
correspondingly,
and let
$$
S(j,l;\nu )=\cases {(-1)^\nu & for $j=l$\cr \hskip 0.3cm 1 & for
$j<l$\cr
-1 &for $j>l$\cr},\quad
\theta(i)=\cases {1 & for
$i\geq 0$ \cr  0 & for $i<0$\cr}, \quad
L_{ij}=M_{ij}-i. \eqno(16)
$$
Set moreover
$$
e^0_i=f_i, \quad e^1_i=e_i, \quad i\in \Z. \e(17)
$$

Let $\{\rho(e_i), \rho(f_i), \rho(h_i), \rho(c) \}_{i\in \Z}$ be a
collection of $\C[[h]]$-endomorphisms of $\V$, defined on any
$C$-pattern $(M)\in \G$, as follows 
(see also (50), (51), (54)
and (55)):

$$
\eqalignno{
& \rho(e_{-1}^{1-\mu})(M)=( [L_{-1,2}-L_{0,1}-\mu]
[L_{0,1}-L_{0,2}+\mu ])^{1/2}(M)_{-(-1)^\mu \{0,1\} }, \;\quad
\mu =0,1\;, & (18)\cr
& &\cr
& \rho(e^{\;\mu}_{(-1)^\nu i-1})(M)=-\sum_{j=1-i-\nu }^{i-1}
\sum_{l=-i}^{i+\nu-1} S(j,l;\nu )& \cr
& &\cr
& \times \Biggl(-{\prod_{k\not= l=-i}^{i+\nu -1}[L_{k,2i+\nu }-
L_{j,2i+\nu -1}
-(-1)^\nu \mu]\prod_{k=1-i}^{i+\nu-2}[L_{k,2i+\nu-2}-L_{j,2i+\nu
-1}-(-1)^\nu
\mu]\over
{\prod_{k\not= j=1-i-\nu}^{i-1}[L_{k,2i+\nu -1}-L_{j,2i+\nu -1}]
[L_{k,2i+\nu -1}-L_{j,2i+\nu-1}+(-1)^{\mu +\nu}]}} & \cr
& &\cr
& \times {\prod_{k=-i-\nu}^{i}[L_{k,2i+\nu +1}-L_{l,2i+\nu }
+(-1)^\nu (1-\mu )]
\prod_{k\not= j=1-i-\nu}^{i-1}[L_{k,2i+\nu-1}-L_{l,2i+\nu }+(-1)^\nu
(1-\mu )]\over
{\prod_{k\not= l=-i}^{i+\nu -1}[L_{k,2i+\nu }-L_{l,2i+\nu }]
[L_{k,2i+\nu }-L_{l,2i+\nu}+(-1)^{\mu +\nu}]}}\Biggr)^{1/2}& \cr
&& \cr
& \times (M)_{-(-1)^{\mu +\nu }\{ j,2i-1+\nu\} }^{-(-1)^{\mu +\nu }\{
l,2i+\nu\} },\quad  i \in {\bf N},\quad  \mu , \nu =0,1\;,&(19) \cr
&& \cr
& \rho(h_{i})(M)=\left(\sum_{j=-|i|}^{|i|+\theta
(i)-1}M_{j,2|i|+\theta(i)} -
\sum_{j=-|i|+1-\theta (i)}^{|i|-1}M_{j,2|i|+\theta
(i)-1}+(\xi_1-\xi_0)\theta (-i)-\xi_1\right) (M), \quad i\in \Z, & (20)\cr
&& \cr
& \rho (c)(M)=(\xi _0-\xi _1)(M). & (21) \cr
}
$$

\bigskip
\noindent
Above and throughout
$
[x]={q^x-q^{-x}\over {q-q^{-1}}}\in \C[[h]].
$
If a pattern from the right hand side of (19)  does not belong to
$\G$, i.e., it is not a $C$-pattern, then the corresponding term
has to be deleted. (The coefficients in front of all such
patterns are  undefined, they contain zero multiples in the
denominators. Therefore an equivalent statement is that all terms
with zeros in the denominators have to be removed).  With this
convention all coefficients in front of the $C-$patterns in r.h.s
of (18)-(21) are well defined as elements from $\C[[h]]$.

\bigskip
\noindent
{\bf Proposition 1.} {\it The endomorphisms $\{\rho(e_i), \rho(f_i),
\rho(h_i),\rho(c)  \}_{i\in \Z}$ satisfy Eqs. (1)-(3) with $\rho(e_i),
\rho(f_i), \rho(h_i)$ and  $\rho(c)$ substituted for 
$e_i, f_i, h_i,$ and $c $, respectively.}

\smallskip\n
{\it Proof.}
The proof is based on a direct verification of the relations
(1)-(3). This verification is lengthy and nontrivial. The most
difficult to check are the last Cartan relations in (1),
corresponding to $i=j$. In order to show they hold, one has to
prove as an intermediate step that the following identities are
fulfilled:

$$
\eqalignno{
&\sum_{s=0}^1 \;
\sum_{j=1-k}^{k-1}\; \sum_{l=-k}^{k-1} (-1)^s
{\prod_{i\neq l=-k}^{k-1}[L_{i,2k}-L_{j,2k-1}+s-1]
\prod_{i=1-k}^{k-2}[L_{i,2k-2}-L_{j,2k-1}+s-1]
\over {\prod_{i\neq j=1-k}^{k-1}[L_{i,2k-1}-L_{j,2k-1}+s]
[L_{i,2k-1}-L_{j,2k-1}+s-1]  }}\times & \cr
&&\cr
& \hskip 22mm  \times {\prod_{i=-k}^k[L_{i,2k+1}-L_{l,2k}+s]
\prod_{i\neq j=1-k}^{k-1}[L_{i,2k-1}-L_{l,2k}+s]
\over{\prod_{i\neq l=-k}^{k-1}[L_{i,2k}-L_{l,2k}+s]
[L_{i,2k}-L_{l,2k}+s-1]}} &\cr
&&\cr
& \hskip 22mm =\left[\sum_{j=-k+1}^{k-1}L_{j,2k-1}-
\sum_{j=-k+1}^{k-2}L_{j,2k-2}-
\sum_{j=-k}^k L_{j,2k+1}+\sum_{j=-k}^{k-1} L_{j,2k}-1\right],
\quad \forall k\in \N, & (22)\cr
}
$$

\bigskip
$$
\eqalignno{
&\sum_{s=0}^1 \;
\sum_{j=-k}^{k-1}\; \sum_{l=-k}^{k} (-1)^s
{\prod_{i\neq l=-k}^{k}[L_{i,2k+1}-L_{j,2k}-s+1]
\prod_{i=1-k}^{k-1}[L_{i,2k-1}-L_{j,2k}-s+1]
\over {\prod_{i\neq j=-k}^{k-1}[L_{i,2k}-L_{j,2k}-s]
[L_{i,2k}-L_{j,2k}-s+1]  }}\times & \cr
&&\cr
& \hskip 20mm  \times {\prod_{i=-k-1}^k[L_{i,2k+2}-L_{l,2k+1}-s]
\prod_{i\neq j=-k}^{k-1}[L_{i,2k}-L_{l,2k+1}-s]
\over{\prod_{i\neq l=-k}^{k}[L_{i,2k+1}-L_{l,2k+1}-s]
[L_{i,2k+1}-L_{l,2k+1}-s+1]}} &\cr
&&\cr
& \hskip 20mm =\left[\sum_{j=-k-1}^{k}L_{j,2k+2}-
\sum_{j=-k}^{k}L_{j,2k+1}-
\sum_{j=-k}^{k-1} L_{j,2k}+\sum_{j=-k+1}^{k-1} L_{j,2k-1}-1\right],
\quad \forall k\in \N. &
(23)
\cr
}
$$

We proceed to verify (22) and (23). Let
$$
\eqalignno{
&F_k(h)=\sum_{s=0}^1 \;
\sum_{j=1-k}^{k-1}\; \sum_{l=-k}^{k-1} (-1)^s
{\prod_{i\neq l=-k}^{k-1}[L_{i,2k}-L_{j,2k-1}+s-1]
\prod_{i=1-k}^{k-2}[L_{i,2k-2}-L_{j,2k-1}+s-1]
\over {\prod_{i\neq j=1-k}^{k-1}[L_{i,2k-1}-L_{j,2k-1}+s]
[L_{i,2k-1}-L_{j,2k-1}+s-1]  }}\times & \cr
&&\cr
& \hskip 22mm  \times {\prod_{i=-k}^k[L_{i,2k+1}-L_{l,2k}+s]
\prod_{i\neq j=1-k}^{k-1}[L_{i,2k-1}-L_{l,2k}+s]
\over{\prod_{i\neq l=-k}^{k-1}[L_{i,2k}-L_{l,2k}+s]
[L_{i,2k}-L_{l,2k}+s-1]}}, &\cr
&&\cr
&G_k(h)=\left[\sum_{j=-k+1}^{k-1}L_{j,2k-1}-
\sum_{j=-k+1}^{k-2}L_{j,2k-2}-
\sum_{j=-k}^k L_{j,2k+1}+\sum_{j=-k}^{k-1} L_{j,2k}-1\right]. &\cr
}
$$
$F_k(h)$ and $G_k(h)$ are formal power series,
$$
F_k(h)=\sum_{i=0}^\infty c_{ki} h^i,\q
G_k(h)=\sum_{i=0}^\infty d_{ki} h^i.
$$
In order to prove that (22) holds one has to show that
$c_{ki}=d_{ki}$ for any $k$ and $i$.
To this end replace $h$ in (22) by a complex variable $x\in \C$,
so that $x$ takes values on any curve $\gamma \subset \C$ and
$x\notin i\pi \Q$. Then $F_k(x)$ and $G_k(x)$ are well defined
for any $x\in \gamma$ and
$$
F_k(x)=\sum_{i=0}^\infty c_{ki} x^i,\q
G_k(x)=\sum_{i=0}^\infty d_{ki} x^i.
$$
Therefore $c_{ki}=d_{ki}$ if $F_k(x)=G_k(x)$.
Since $h$ appears in (22) only through $q=e^{h/2}$,
we have to show that (22) is an identity  for $q$ being a
number, which is not a root of 1. The same arguments hold also
for (23).

Let $q\in \C$ be not a root of 1. Setting
$q^{2L_{i,2k-1}}=A_i, \;\; q^{2L_{i,2k}}=B_i, \;\;
q^{2L_{i,2k+1}}=C_i,
\;\; q^{2L_{i,2k-2}}=D_i, \;\; 2k=n$ in (22) and
$q^{2(L_{i,2k}-1)}=A_i, \;\; q^{2L_{i,2k+1}}=B_i, \;\;
q^{2(L_{i,2k+2}-1)}
=C_i,
\;\; q^{2L_{i,2k-1}}=D_i, \;\; 2k+1=n$ in (23)
and relabeling the indices in an appropriate way, one reduces both
identities to the following one:

\bigskip

$$ \sum_{j=1}^{n-1} \sum_{l=1}^n
{q \prod_{i\neq l=1}^n (A_j-q^{-2}B_i)
\prod_{i=1}^{n-2}(A_j-q^{-2}D_i)
\prod_{i=1}^{n+1}(B_l-C_i)
\prod_{i\neq j=1}^{n-1}(B_l-A_i)
\over {A_jB_l\prod_{i\neq j=1}^{n-1}(A_j-A_i)(A_j-q^{-2}A_i)
\prod_{i\neq l=1}^n (B_l-B_i)(B_l-q^{-2}B_i)}}$$

$$ -\sum_{j=1}^{n-1} \sum_{l=1}^n
{q^{-1} \prod_{i\neq l=1}^n (A_j-B_i)
\prod_{i=1}^{n-2}(A_j-D_i)
\prod_{i=1}^{n+1}(B_l-q^2C_i)
\prod_{i\neq j=1}^{n-1}(B_l-q^2A_i)
\over {A_jB_l\prod_{i\neq j=1}^{n-1}(A_j-A_i)(A_j-q^{2}A_i)
\prod_{i\neq l=1}^n (B_l-B_i)(B_l-q^{2}B_i)}}$$

$$
=(q-q^{-1})\left( 1-q^2{\prod_{i=1}^{n-2}D_i \prod_{i=1}^{n+1} C_i
\over {\prod_{i=1}^{n-1}A_i \prod_{i=1}^n B_i}}\right), \eqno (24)
$$
which can be written also in the form

$$ \sum_{j=1}^{n-1}
{q \prod_{i=1}^n (A_j-q^{-2}B_i)
\prod_{i=1}^{n-2}(A_j-q^{-2}D_i)\over
{A_j\prod_{i\neq j=1}^{n-1}(A_j-A_i)(A_j-q^{-2}A_i)}}
\sum_{l=1}^n{
\prod_{i=1}^{n+1}(B_l-C_i)
\prod_{i\neq j=1}^{n-1}(B_l-A_i)
\over {(A_j-q^{-2}B_l)B_l
\prod_{i\neq l=1}^n (B_l-B_i)(B_l-q^{-2}B_i)}}$$

$$- \sum_{l=1}^n
{q^{-1}\prod_{i=1}^{n+1}(B_l-q^2C_i)
\prod_{i=1}^{n-1}(B_l-q^2A_i)
\over {B_l
\prod_{i\neq l=1}^n (B_l-B_i)(B_l-q^{2}B_i)}}
\sum_{j=1}^{n-1}
{\prod_{i\neq l=1}^n (A_j-B_i)
\prod_{i=1}^{n-2}(A_j-D_i)
\over {A_j(B_l-q^2A_j)
\prod_{i\neq j=1}^{n-1}(A_j-A_i)(A_j-q^{2}A_i)}}
$$

$$=(q-q^{-1})\left( 1-q^2{\prod_{i=1}^{n-2}D_i \prod_{i=1}^{n+1} C_i
\over {\prod_{i=1}^{n-1}A_i \prod_{i=1}^n B_i}}\right). \eqno (25)$$

\bigskip
\noindent
Consider the complex functions:

$$f_1(z)={\prod_{i=1}^{n+1}(z-C_i)\prod_{i\neq j=1}^{n-1}(z-A_i)
\over {(z-q^2A_j)\prod_{i=1}^n (z-B_i)(z-q^{-2}B_i)}}, \quad
f_2(z)={\prod_{i\neq l=1}^{n}(z-B_i)\prod_{i=1}^{n-2}(z-D_i)
\over {(z-q^{-2}B_l)\prod_{i=1}^{n-1} (z-A_i)(z-q^{2}A_i)}}.
\eqno (26)$$

\bigskip
\noindent
The function $f_1(z)$ (resp. $f_2(z)$) is holomorphic over $\C$
except in
its simple poles
$q^2A_j,\; B_1,\ldots , B_n,$  $q^{-2}B_1, \ldots, q^{-2}B_n$ (resp.
$q^{-2}B_l,\; A_1, \ldots , A_{n-1},\; q^2A_1, \ldots , q^2A_{n-1}.$)
Let $C$ be a closed curve whose interior contains all poles of $f_1(z)$
(resp. $f_2(z).$)
The Residue theorem of complex analysis
implies that $\oint _{\C}f_j(z)dz=
2\pi i\sum Res(f_j(z)),\; j=1,2.$ On the other hand,

$$\oint _{\C}f_j(z)dz=-2\pi i Res f_j(\infty ) \quad j=1,2. \eqno(27)$$

\noindent
Applying the Residue Theorem for both functions (26) and inserting the
results in (25) one obtains:

$$\sum_{j=1}^{n-1}{\prod_{i=1}^{n-2}(A_j-q^{-2}D_i)
\prod_{i=1}^{n+1}(A_j-q^{-2}C_i)
\over {A_j \prod_{i\neq j=1}^{n-1}(A_j-A_i)
\prod_{i=1}^{n}(A_j-q^{-4}B_i)}}
+\sum_{l=1}^n{\prod_{i=1}^{n+1}(B_l-q^{2}C_i)
\prod_{i=1}^{n-2}(B_l-q^{2}D_i)
\over {B_l \prod_{i\neq l=1}^{n}(B_l-B_i)
\prod_{i=1}^{n-1}(B_l-q^{4}A_i)}}
$$

$$= 1-q^2{\prod_{i=1}^{n-2}D_i \prod_{i=1}^{n+1} C_i
\over {\prod_{i=1}^{n-1}A_i \prod_{i=1}^n B_i}}. \eqno(28)$$

\bigskip
\noindent
In order to prove this last identity consider the function:
$$
f(z)={\prod_{i=1}^{n-2}(z-q^{-2}D_i)\prod_{i=1}^{n+1}(z-q^{-2}C_i)
\over {z\prod_{i=1}^{n-1} (z-A_i)\prod_{i=1}^n(z-q^{-4}B_i)}}. \eqno
(29)
$$

\noindent
Applying again the Residue Theorem on this function one obtains
(28) which proves  formulas (22) and (23).

\smallskip
In addition to the identities (22) and (23) the last Cartan
relation (1) is valid if:

$$\sum_{i=1}^n\left( {\prod_{j=1}^{n-1}[a_i-b_j-1]
\prod_{j=1}^{n-1}[a_i-c_j-1]
\over {\prod_{j\neq i=1}^n [a_i-a_j][a_i-a_j-1]}}-
{\prod_{j=1}^{n-1}[a_i-b_j]
\prod_{j=1}^{n-1}[a_i-c_j]
\over {\prod_{j\neq i=1}^n [a_i-a_j][a_i-a_j+1]}}\right) =0, \eqno(30)
$$
which is proved in a similar way. The verification of all other
relations (1)-(3) is based on (30) or/and  on some of the following
identities:
$$ \eqalignno{
&([a-b-1][c-b-1]-[2][a-b][c-b-1]+[a-b][c-b])\left({[a-d][c-e-1]\over
{[d-e-1][c-a-1]}}+ {[c-d-1][a-e]\over {[d-e+1][c-a-1]}}\right)+ & \cr
&& \cr
& ([a-b-1][c-b-1]-[2][a-b-1][c-b]+[a-b][c-b])
\left({[a-e-1][c-d]\over {[d-e-1][c-a+1]}}+
{[a-d-1][c-e]\over {[d-e+1][c-a+1]}}\right)=0, & \cr
&& (31) \cr
&{[a-1][b-1]-[2][a][b-1]+[a][b]\over {[a-b+1]}}+
{[a-1][b-1]-[2][a-1][b]+[a][b]\over {[a-b-1]}}=0, & (32) \cr
&& \cr
& [a-1]-[2][a]+[a+1]=0.  & (33) \cr
}
$$ 
\noindent
This completes the proof. \hfill $[]$

\n
{\it Remark.} We may have constructed the endowmorphisms
$\{\rho(e_i),\rho(f_i),\rho(h_i),\rho(c) \}_{i\in \Z}$ from the
results on the representations of $\gl$, announced in Ref. [14]
and, more precisely, from the endowmorphisms
$\{\rho(E_i),\rho(F_i),\rho(H_i)\}_{i\in \Z}$ of the $\Ch-$module
$V(\{M\})$, which satisfy the Cartan and the Serre relations for
$\gl$ (Ref. [14], Eqs. (16)-(18)). This possibility is based on the
observation that the $C[[h]]-$linear map $\varphi$,
defined on the generators as
$$
\eqalign{
& \varphi (E_i)=e_i, \quad \varphi (F_i)=f_i, \quad
\varphi (c)=c, \cr
& \varphi (H_i)=h_i+(\theta (-i)+\alpha )c, \q \alpha \in \C[h]  \cr
}\e(34)
$$
and extended by associativity is an (algebra) isomorphism of
$U_h(gl_\infty) \oplus \C [[h]] c$ onto $U_h(a_\infty ).$
Then the endomorphisms
$\{\rho(e_i),\rho(f_i),\rho(h_i),\rho(c) \}_{i\in \Z}$
defined according to (34) (with $\alpha (\xi_0-\xi_1)=\xi_1$) as
$$
\eqalign{
& \rho (e_i)=\rho (E_i), \quad \rho (f_i)=\rho (F_i), \quad
\rho (c)=\xi_0-\xi_1, \cr
& \rho (h_i)=\rho (H_i)-((\xi_0-\xi_1)\theta (-i)+\xi_1 ), \cr
}\e(35)
$$
obey Eqs. (18)-(21). Here we have given a direct proof of
{\it Proposition 1}, since the corresponding {\it Proposition 1} in
Ref. [14] was only stated. Its proof would have been based
again on the identities (22) and (23).

Consider $\rho$ as a $\Ch-$linear operator from $\Ua$ into
$End \; \V$. So far $\rho$ is defined only on the Chevalley
generators $\{e_i, f_i, h_i, c  \}_{i\in \Z}$.
Extend the domain of its definition on ${\hat U}_h(a_\infty)$:
if $\rho$ has already been defined on $a, b \in {\hat
U}_h(a_\infty)$, then set
$$
\rho(\alpha a + \beta b)=\alpha \rho(a) + \beta \rho(b),\quad
\rho(ab)=\rho(a)\rho(b), \quad a, b \in {\hat U}_h(a_\infty), \quad
\alpha, \beta \in \C[[h]]. \e(36)
$$
As we know from (10), any element $u\in \Ua$ can be represented
as a
sum $u=\sum_{i=0}^\infty u_i h^i,
\; u_i \in {\hat U}_h(a_\infty)$. Then for an arbitrary
$v\in \V$, writing it as in (14), we have
$$
\left(\sum_{i=0}^\infty \rho(u_i)h^i\right)v=
\left(\sum_{i=0}^\infty \rho(u_i)h^i\right)
\left(\sum_{j=0}^\infty v_jh^j\right)=
\sum_{n=0}^\infty \left(\sum_{m=0}^n \rho(u_{n-m})v_m\right)h^n
\in \V, \eqno(37)
$$
since
$$
\sum_{m=0}^n \rho(u_{n-m})v_m \in {\hat V}(\{M\};\xi_0,\xi_1).
$$
Using (37), we extend $\rho$ on $U_h(a_\infty)$:
$$
\rho(u)=\sum_{i=0}^\infty \rho(u_i)h^i\in End \; \V
\;\;\; \forall \;\; u\in U_h(a_\infty).
\e(38)
$$
Hence $\rho$ is a well defined map from $\Ua$ into
$End\;\V$,
$$
\rho: \; \Ua \rightarrow End\;\V. \e(39)
$$

\bigskip
\noindent
{\bf Proposition 2.} {\it The map (39), acting on the $C-$basis according
to Eqs. (18)-(21), defines a highest weight irreducible representation of 
$\Ua$ in $\V)$.}

\smallskip\noindent
{\it Proof.} According to {\it Proposition 1}, 
(36), (38), $\rho$ is a
$\C[[h]]$-homomorphism of $\Ua$ in 
$End\;\V$. It is continuous in
the $h-$adic topology. Indeed, let $u\in \Ua$. Then any
neighbourhood $W(\rho(u))$ of $\rho(u)$ contains a basic
neighbourhood $\rho(u)+h^nEnd\;\V \subset W(\rho(u)) $.
Evidently
$$
\rho\left( u+h^n \Ua \right) \subset
\rho(u)+h^n End\;\V \subset W(\rho(u))
\e(40)
$$
and therefore $\rho$ is continuous in $u$ for any $u\in \Ua$.
Hence $\V$ is a $\Ua-$module. It is a highest weight module with
respect to the "Borel" subalgebra $U_h(n_+)$. The highest weight
vector $(\hat{M})$, which by definition satisfies the condition
$\rho(U_h(n_+))(\hat{M})=0$ and is an eigenvector of $\rho(H)$,
corresponds to the one from (11) with
$$
\hat{M}_{i, 2k+\theta -1}=M_i , \quad \forall \; k \in \N, \quad \theta
=0,1,
\quad
\;\; i \in [-\theta -k+1, k-1]. \eqno(41)
$$
Moreover, $\V=\rho(\Ua)(\hat{M})$.  Since
$\V$ contains no other singular vectors, vectors
annihilated from $\rho(U_h(n_+))$, each $\V$ is an irreducible
$\Ua$-module. The proof of the latter follows from the
results in Ref. [12] and the observation that each (deformed)
matrix element in the transformation relations (18)-(21) is zero
only if the corresponding nondeformed matrix element vanishes.
\hfill $[]$
\vskip 24pt

\n
{\bf 3. REPRESENTATIONS OF THE ALGEBRA U$_h$($A_\infty$)}

\bigskip\n
In this section we show that some of the $\Ua-$modules $\V$ can
be turned into irreducible highest weight $\UA-$modules. First,
following Ref. [8], we recall the definition of the quantum
algebra $\UA$ (denoted by the authors as $U_h(g'(A_\infty))$)
only within the algebra sector (for the other Hopf algebra maps
see Ref. [8]). $\UA$  consists of all elements
$$
u=\sum_{k=0}^\infty h^k \sum_{l=0}^{ l(k) } c^l
\sum_{\alpha , \beta \in Q_+'} \; \sum _{\gamma (k,l)\in \Z_+^\infty }
\sum_{t=1}^{t(\alpha,\beta)}
{\cal F}_{\alpha ,k,l,t}\prod_{i\in \Z} h_{i}^{\gamma (k,l)_i}
{\cal E} _{\beta ,k,l,t}, \;\; infinite \; sums \; over
\; \alpha,\; \beta, \; \gamma,   \e(42)
$$
however certain conditions on the pairs $(\a, \g)$ corresponding
to the non-zero summands are imposed. In order to state them, set
for $\a=\sum_i m_i \a_i \in \Q_+',\; \g=\{\g_i\}_{i\in \Z} 
\in \Z_+^\infty,$
$$
S(\a)=\{i\;|\; m_i\ne 0\},\q S(\g)=\{i\;|\;\g_i\ne 0\},\q
S(\a,\g)=S(\a)\cup S(\g).\e(43)
$$
Connecting $i$ and $j$ if $|i-j|=1$, one can view $S(\a,\g)$ as a
graph. Denote by ${\cal F}(\a,\g)$ the collection of its
connected components. For any $u=\sum_{k=0}^\infty h^k u_k $ as given in (42)
consider the series
$$
u_k= \sum_{l=0}^{ l(k) } c^l
\sum_{\alpha , \beta \in Q_+'} \; \sum _{\gamma (k,l)\in \Z_+^\infty }
\sum_{t=1}^{t(\alpha,\beta)}
{\cal F}_{\alpha ,k,l,t}\prod_{i\in \Z} h_{i}^{\gamma (k,l)_i}
{\cal E} _{\beta ,k,l,t} , \e(44)
$$
and let
$$
{\cal F}(u,k)=\cup {\cal F}(\a,\g), \e(45)
$$
where the union is taken over all $\a$ and $\g$, which appear in
the non-zero summands of (44). For $i\in \Z$ and $k
\in \Z_+$ set
$$
Int(u,k,i)=\{I\in {\cal F}(u,k)|\;i\in I\}.\e(46)
$$

For $r\in \N$ define the series $u(r)$, corresponding to $u$, by
substituting 0 for all $h_i$ ($i\le -r$ or $i>r$) and
for all $e_i,\; f_i$ ($|i|\ge r$).

\bigskip\n
{\bf Definition 1.} $[8]$ {\it The series $u$ of the form (42) is said
to belong to $\UA$, provided

($i$) for any $k\in \Z_+$ and any $i\in \Z$ the set $Int(u,k,i)$
   is finite;

($ii$) $u(r)\in \Ua$ for all $r \in \N$.}

\bigskip

We turn to construct a class of representations of $\UA$. The idea
is to show that within certain $\V$ the domain of the definition  of
the operator $\rho: \; \Ua \rightarrow End\;\V$  (see (39)) can
be extended from $\Ua$ to $\UA$, so that $\rho$ is a
representation of $\UA$ in $\V$.  The construction is a natural
one. We assume that $\rho (u)$ is a continuous $\Ch-$linear
operator in $\V$ for any $u \in \UA$.  Then for any
$v=\sum_{j=0}^\infty v_jh^j\in\V,\;v_j \in {\hat
V}(\{M\};\xi_0,\xi_1)$
$$
\rho(u)v=\sum_{j=0}^\infty \left(\rho(u)v_j\right)h^j. \e(47)
$$
The above series is well defined, if $\rho(u)v_j \in \V$.
Since $v_j$ is a (finite) $\Ch-$linear combination of $C-$basis
vectors the latter holds, and hence $\rho(u)$ is defined as an
operator in $\V$, if $\rho(u)(M)\in \V$ for any $C-$pattern
$(M)$.  Thus, the first step is to clarify which are the
$\Ua-$modules $\V$, for which this can be achieved.

\bigskip\n
{\bf Definition 2.} {\it We say that $\V$ is of a finite-signature or,
more precisely, of $(m,n)$ signature and write
$\{M\}=\{M\}_{m,n}$ if there exist integers $m\le n
\in \Z$, such that $M_m=M_{m-k}$ and $M_n=M_{n+k}$ for all
$k\in \N$.}

\bigskip
We now proceed to show that each finite signature $\Ua$ module
$\Vf$ can be turned into a $\UA$ module. To this end we prove first
a few preliminary propositions.

Denote by $\V_N, \; 1<N\in \N$ the subspace of $\V$,
which is a $\Ch-$linear envelope of all $C-$basis vectors $(M)$,
for which
$$
M_{i,2k+\t-1}=M_i \;\; \forall \; 2k+\t-1 \geq N, \quad
 \theta =0,1, \;\; i \in [1-\theta -k, k-1]; \e(48)
$$
Note that $\V_N$ is a finite dimensional subspace.

\smallskip\n
{\bf Proposition 3.}
$$
\rho(e_k)\V_N=0, \;\; if \;\; k \notin
(-{1\over 2}(N+1), {1\over 2}(N-2)).
\e(49)
$$

\smallskip\n
{\it Proof.} From (19) one obtains:

$$
\eqalignno{
\rho(e_k)(M)=&-\sum_{j=-k}^{k} \; \sum_{l=-k-1}^{k} S(j,l;0 )&\cr
&&\cr
& \times \left|
{\prod_{i\not= l=-k-1}^{k}[L_{i,2k+2}-L_{j,2k+1}-1]
\prod_{i=-k}^{k-1}[L_{i,2k}-L_{j,2k+1}-1]\over
{\prod_{i\not= j=-k}^{k}[L_{i,2k+1}-L_{j,2k+1}]
[L_{i,2k+1}-L_{j,2k+1}-1]}}\right|^{1/2} & (50)\cr
&&\cr
& \times\left|{\prod_{i=-k-1}^{k+1}[L_{i,2k+3}-L_{l,2k+2}]
\prod_{i\not= j=-k}^{k}[L_{i,2k+1}-L_{l,2k+2}]\over
{\prod_{i\not= l=-k-1}^{k}[L_{i,2k+2}-L_{l,2k+2}]
[L_{i,2k+2}-L_{l,2k+2}-1]}}\right|^{1/2} (M)_{\{ j,2k+1\} }^{\{
l,2k+2\} }, \; k\ge 0.& \cr
&&\cr
&&\cr
\rho(e_{-k})(M)=&-\sum_{j=-k+1}^{k-2} \; \sum_{l=-k+1}^{k-1} S(j,l;1)&\cr
&&\cr
&\times \left|{\prod_{i\not= l=-k+1}^{k-1}[L_{i,2k-1}-L_{j,2k-2}+1]
\prod_{i=2-k}^{k-2}[L_{i,2k-3}-L_{j,2k-2}+1]\over
{\prod_{i\not= j=-k+1}^{k-2}[L_{i,2k-2}-L_{j,2k-2}]
[L_{i,2k-2}-L_{j,2k-2}+1]}}\right|^{1/2} & (51)\cr
&&\cr
& \times \left|{\prod_{i=-k}^{k-1}[L_{i,2k}-L_{l,2k-1}]
\prod_{i\not= j=-k+1}^{k-2}[L_{i,2k-2}-L_{l,2k-1}]\over
{\prod_{i\not= l=-k+1}^{k-1}[L_{i,2k-1}-L_{l,2k-1}]
[L_{i,2k-1}-L_{l,2k-1}+1]}}\right|^{1/2} (M)_{-\{ j,2k-2\} }^{-\{
l,2k-1\} },\; k>1& \cr
}
$$
If $(M)\in \V_N$ and
$$
\eqalign
{
& if \q k\ge {1\over 2}(N-2)),\;\;then \;\; L_{i,2k+3}=L_{i,2k+2}=L_i=M_i-i,\cr
& if \q 
k\le -{1\over 2}(N+1),\;\; then \;\;L_{i,2k-1}=L_{i,2k}=L_i=M_i-i.\cr
}\e(52)
$$
In both cases the r.h.s. of (50) and (51) contain zero multiples
$(L_l-L_l)$ and therefore vanish. 
\hfill $[]$
\smallskip\n
{\bf Proposition 4.}
$$
\rho(f_k)\Vf_N=0, \;\; if \;\; k \notin
(min\{-{1\over 2} (N+3),m-1 \},max\{{1\over 2} N, n \}).\e(53)
$$

\smallskip\n
{\it Proof.} Let $(M)\in \Vf_N$.
The relations that follow from (19) in this case are
$$
\eqalignno{
\rho(f_k)(M)=&-\sum_{j=-k}^{k}\sum_{l=-k-1}^{k} S(j,l;0 )& \cr
&&\cr
& \times \left|{\prod_{i\not= l=-k-1}^{k}[L_{i,2k+2}-L_{j,2k+1}]
\prod_{i=-k}^{k-1}[L_{i,2k}-L_{j,2k+1}]\over
{\prod_{i\not= j=-k}^{k}[L_{i,2k+1}-L_{j,2k+1}]
[L_{i,2k+1}-L_{j,2k+1}+1]}}\right|^{1/2}& (54) \cr
&&\cr
& \times \left|{\prod_{i=-k-1}^{k+1}[L_{i,2k+3}-L_{l,2k+2}+1]
\prod_{i\not= j=-k}^{k}[L_{i,2k+1}-L_{l,2k+2}+1]\over
{\prod_{i\not= l=-k-1}^{k}[L_{i,2k+2}-L_{l,2k+2}]
[L_{i,2k+2}-L_{l,2k+2}+1]}}\right|^{1/2} (M)_{-\{ j,2k+1\} }^{-\{
l,2k+2\} }, \; k\in \Z_+,& \cr
&&\cr
&&\cr
\rho(f_{-k})(M)=&-\sum_{j=-k+1}^{k-2}\sum_{l=-k+1}^{k-1} S(j,l;1 )&\cr
&&\cr
&\times\left|{\prod_{i\not= l=-k+1}^{k-1}[L_{i,2k-1}-L_{j,2k-2}]
\prod_{i=2-k}^{k-2}[L_{i,2k-3}-L_{j,2k-2}]\over
{\prod_{i\not= j=-k+1}^{k-2}[L_{i,2k-2}-L_{j,2k-2}]
[L_{i,2k-2}-L_{j,2k-2}-1]}}\right|^{1/2} &(55) \cr
&&\cr
& \times \left|{\prod_{i=-k}^{k-1}[L_{i,2k}-L_{l,2k-1}-1]
\prod_{i\not= j=-k+1}^{k-2}[L_{i,2k-2}-L_{l,2k-1}-1]\over
{\prod_{i\not= l=-k+1}^{k-1}[L_{i,2k-1}-L_{l,2k-1}]
[L_{i,2k-1}-L_{l,2k-1}-1]}}\right|^{1/2} (M)_{\{ j,2k-2\} }^{\{
l,2k-1\}}, \; k>1. &\cr
}
$$
If $k\ge N/2$, then
$L_{i,2k+3}=L_{i,2k+2}=L_{i,2k+1}=L_{i,2k}=L_{i} 
=M_i-i$ and therefore (54)
reads:
$$
\eqalign{
\rho(f_k)(M)=& -\sum_{j=-k}^{k}\sum_{l=-k-1}^{k} S(j,l;0 )
\left|{\prod_{i\not= l=-k-1}^{k}[L_{i}-L_{j}]
\prod_{i=-k}^{k-1}[L_{i}-L_{j}]\over
{\prod_{i\not= j=-k}^{k}[L_{i}-L_{j}]
[L_{i}-L_{j}+1]}}\right|^{1/2} \cr
&\cr
& \times \left|{\prod_{i=-k-1}^{k+1}[L_{i}-L_{l}+1]
\prod_{i\not= j=-k}^{k}[L_{i}-L_{l}+1]\over
{\prod_{i\not= l=-k-1}^{k}[L_{i}-L_{l}]
[L_{i}-L_{l}+1]}}\right|^{1/2} (M)_{-\{ j,2k+1\} }^{-\{
l,2k+2\} } \cr
}
$$
Above only the term with $j=k$ survives:
$$
\eqalign{
\rho(f_k)(M)=&
-\sum_{l=-k-1}^{k} S(k,l;0 )
\left|{\prod_{i\not= l=-k-1}^{k}[L_{i}-L_{k}]
\prod_{i=-k}^{k-1}[L_{i}-L_{k}]\over
{\prod_{i=-k}^{k-1}[L_{i}-L_{k}]
[L_{i}-L_{k}+1]}}\right|^{1/2} \cr
&\cr
& \times \left|{\prod_{i=-k-1}^{k+1}[L_{i}-L_{l}+1]
\prod_{i=-k}^{k-1}[L_{i}-L_{l}+1]\over
{\prod_{i\not= l=-k-1}^{k}[L_{i}-L_{l}]
[L_{i}-L_{l}+1]}}\right|^{1/2} (M)_{-\{ k,2k+1\} }^{-\{
l,2k+2\} } \cr
&\cr
=-S(k,k;0 )&
\left|{\prod_{i=-k-1}^{k-1}[L_{i}-L_{k}]
\prod_{i=-k-1}^{k+1}[L_{i}-L_{k}+1]
\prod_{i=-k}^{k-1}[L_{i}-L_{k}+1]\over
{\prod_{i=-k}^{k-1}
[L_{i}-L_{k}+1]
\prod_{i=-k-1}^{k-1}[L_{i}-L_{k}]
[L_{i}-L_{k}+1]}}\right|^{1/2} (M)_{-\{ k,2k+1\} }^{-\{
k,2k+2\} } \cr
&\cr
-\sum_{l=-k-1}^{k-1} S(k,l;0 )&
\left|{\prod_{i\not= l=-k-1}^{k}[L_{i}-L_{k}]
\prod_{i=-k-1}^{k+1}[L_{i}-L_{l}+1]
\prod_{i=-k}^{k-1}[L_{i}-L_{l}+1]\over
{\prod_{i=-k}^{k-1}
[L_{i}-L_{k}+1]
\prod_{i\not= l=-k-1}^{k}[L_{i}-L_{l}]
[L_{i}-L_{l}+1]}}\right|^{1/2} (M)_{-\{ k,2k+1\} }^{-\{
l,2k+2\} } . \cr
}
$$
In the last term $l\ne k$. Hence
$$
\prod_{i\not= l=-k-1}^{k}[L_{i}-L_{k}] =[L_k-L_k]
\prod_{i\not= l=-k-1}^{k-1}[L_{i}-L_{k}]=0
$$
and therefore it vanishes. Then
$$
\rho(f_k)(M)=-
\left|{
\prod_{i=-k-1}^{k+1}[L_{i}-L_{k}+1]
\over {\prod_{i=-k-1}^{k-1}
[L_{i}-L_{k}+1]}}\right|^{1/2} (M)_{-\{ k,2k+1\} }^
{-\{ k,2k+2\} }=
-\left|{ [L_{k+1}-L_{k}+1]}\right|^{1/2} (M)_{-\{ k,2k+1\} }^{-\{
k,2k+2\} },
$$
i.e.,
$$
\rho(f_k)(M)=-\left|{ [M_{k+1}-M_{k}]}\right|^{1/2}
(M)_{-\{ k,2k+1\} }^{-\{k,2k+2\} }.\e(56)
$$
If $k\geq n$ then $M_{k+1}=M_k$ and
$\rho (f_k)(M)=0$.
In a similar way one derives from (55) that
$\rho (f_k)(M)=0$ if
$k \le min\{-{1\over 2} (N+3),m-1 \}$, which proves (53).
\hfill $[]$
\smallskip\n
{\bf Proposition 5.}
$$
\rho(h_k)\Vf_N=0, \;\; if \;\; k \notin
(min\{-{1\over 2} (N+1),m \},max\{{1\over 2} N, n \}).\e(57)
$$

\smallskip
The proof follows easily from (20).

Set
$$
r_N=max\{{1\over 2} (N+3),1-m, n \}.\e(58)
$$
\n
From the last three propositions one concludes.

\smallskip\n
{\bf Corollary 1.} {\it If} $k \notin (-r_N,r_N)$, {\it then}
$$
\rho(h_k)\Vf_N=\rho(e_k)\Vf_N=\rho(f_k)\Vf_N=0.\e(59)
$$

\smallskip\n
{\bf Proposition 6.} 
$\rho(U_h(n_+)_{\b})\V_N \subset \V_N$ {\it for any} $N\in \N$.
{\it More precisely,}
$$
\eqalignno{
& If \;\; S(\b) \subset (-{(N+1)\over 2}, {(N-2)\over
2})\equiv I_N, \;\; then \;\; \rho(U_h(n_+)_{+\b})\V_N \subset \V_N. &(60)\cr
&&\cr
&If \;\; S(\b) \not\subset \; I_N, \;\;
then \;\; \rho(U_h(n_+)_{\b}) \V_N=0. & (61) \cr 
}
$$

\smallskip\n
{\it Proof.} From (19) (or directly from (50) and (51)) one
concludes that
$$
\rho(e_j)(M)\in \V_N ,\;\; \forall \;\; j \in
(-{(N+1)\over 2}, {(N-2)\over 2})\;\;and \;\; \forall \;\; (M)\in \V_N.
\e(62)
$$
Hence (60) holds.

Assume $S(\b) \not\subset \; I_N$ and
let $P_\b$ be a monomial of the generators $e_i, \;
i\in S(\b)$. $P_\b$ can be represented as $P_\b=Q'e_n Q$,
where $Q$ depends only on the generators $e_j$ with $j \in
S(\b)\cap I_N$ and $n \notin I_N$. Then (62) yields that
$\rho(Q)(M)\in \V_N$  and therefore ({\it Proposition 3}) 
$\rho(P_\b)(M)=\rho(Q')\rho(e_n)\rho(Q)(M)=0$.
\hfill $[]$

\smallskip\n
{\bf Proposition 7.} {\it Let $u$ be any element from $\UA$,
represented as in (42). Then} 
$$
\sum_{k=0}^\infty h^k\; \sum_{l=0}^{ l(k) } \rho(c)^l
\sum_{\alpha , \beta \in Q_+'} \; \sum _{\gamma (k,l)\in \Z_+^\infty }
\sum_{t=1}^{t(\alpha,\beta)}
\rho({\cal F}_{\alpha ,k,l,t})\prod_{i\in \Z} \rho(h_{i})^{\gamma (k,l)_i}
\rho({\cal E} _{\beta ,k,l,t})(M)\in \Vf \e(63)
$$
{\it for any 
$(M)$ from the basis $ \G$ of $\Vf $. 
For a fixed $k$ the number of the non-zero
summands in (63) is finite. }

\smallskip\n
{\it Proof.} In order to prove the proposition it 
suffices to
show that 
$$
\sum_{\alpha , \beta \in Q_+'} \; \sum _{\gamma (k,l)\in \Z_+^\infty }
\sum_{t=1}^{t(\alpha,\beta)}
\rho({\cal F}_{\alpha ,k,l,t})\prod_{i\in \Z} \rho(h_{i})^{\gamma (k,l)_i}
\rho({\cal E} _{\beta ,k,l,t})(M)\in \Vh .\e(64)
$$
Let for definiteness $(M)\in \Vf_N$ and let $\b_0$ be
the weight of (M). Since ${\cal E} _{\beta ,k,l,t}\in
U_h(n_+)_{\b}$, $\;\rho({\cal E} _{\beta ,k,l,t})(M)\in \Vf_N$ ({\it
Proposition 6}). Each $\rho({\cal E} _{\beta ,k,l,t})(M)$ has a weight
$\b+\b_0$ and the nonzero vectors $\rho({\cal E} _{\beta ,k,l,t})(M)$,
corresponding to different $\b \in \Q_+'$, are linearly independent
(over \Ch).  Since $\Vf_N$ is a finite dimensional subspace,
$\rho({\cal E} _{\beta ,k,l,t})(M)\ne 0 $ only for a finite number of
$\b \in \Q_+'$. Hence the sum over $\b$ in (64) is finite.
Setting $\rho({\cal E} _{\beta ,k,l,t})(M)=v_{\b,k,l,t} \in
\Vh_N$, we obtain for the l.h.s. of (64)
$$
\sum_{\alpha , \beta \in Q_+'} \; \sum _{\gamma (k,l)\in \Z_+^\infty }
\sum_{t=1}^{t(\alpha,\beta)}
\rho({\cal F}_{\alpha ,k,l,t})\prod_{i\in \Z} \rho(h_{i})^{\gamma (k,l)_i}
v_{\b,k,l,t}. \e(65)
$$

We proceed to show that  the sum over $\a$ and $\g$ in (65) is finite
too. Without loss of generality we assume that every
${\cal F}_{\alpha ,k,l,t}\prod_{i\in \Z}h_i^{\gamma (k,l)_i}$ 
is a monomial of $\{f_i, h_i\}_{i\in
\Z}$. Consider any term from (65), corresponding to a particular
pair $(\a, \g)$:
$$
\rho({\cal F}_{\alpha ,k,l,t})\prod_{i\in \Z_+} \rho(h_{i})^{\gamma (k,l)_i}
v_{\b,k,l,t}. \e(66)
$$
Let ${\cal F}(\a,\g)$ be the connected components of $S(\a,\g)$:

$$
{\cal F}(\a,\g)=\{I_{a_i,b_i} \equiv [a_i,b_i] \;|\;i=1,2,\l,n \},\q
a_i\le b_i\in \Z,\;\; |a_i-b_j|>1, \;if \; i\ne j. \e(67)
$$
In (67) we do not distinguish between a connected component
$I_{a_i,b_i}$ (with a beginning in $a_i$ and end in $b_i$) and
the corresponding to it finite integer interval $[a_i,b_i]$.
Then
$$
{\cal F}_{\alpha ,k,l,t}\prod_{i\in \Z_+} (h_{i})^{\gamma (k,l)_i}=
\prod_{i=1}^n {\cal F}_i,\e(68)
$$
where ${\cal F}_i$ is a monomial of $f_j,h_j,\;j\in [a_i,b_i] $. 
From (1) and (3) it follows that the multiples ${\cal F}_i$
in (68) commute. If $[a_i,b_i] \subset (-\infty, -r_N]$ or
$[a_i,b_i] \subset [r_N, \infty)$, then {\it Corollary 1} yields
that $\rho({\cal F}_i)v_{\b,k,l,t}$=0. Hence also
$
\rho({\cal F}_{\alpha ,k,l,t})\prod_{i\in \Z_+} \rho(h_{i})^{\gamma (k,l)_i}
v_{\b,k,l,t}=0.
$
Thus, the sum in (64) is over such pairs $(\a,\g)$, for which all
connected components of $S(\a,\g)$, namely the elements from
${\cal F}(\a,\g)$, have nonzero intersection with $(-r_N,r_N)$.
There is only a finite number of pairs $(\a,\g)$ with this property,
for which $Int(u,k,i)$ is finite.

The conclusion is that the l.h.s. of the series (64) contains a
finite number of non-zero summands. Since the generators of
$\UA$ (see (18)-(21)) transform $\Vh$ into itself, (64) holds.
Hence (63) holds too. 
\hfill $[]$

Based on {\it Proposition 7}, we extend the domain of the operator
$\rho$ on $\UA$, setting for any $u\in \UA$ (see (42))
$$
\rho(u)=\sum_{k=0}^\infty h^k \sum_{l=0}^{ l(k) } \rho(c)^l
\sum_{\alpha , \beta \in Q_+'} \; \sum _{\gamma (k,l)\in \Z_+^\infty }
\sum_{t=1}^{t(\alpha,\beta)}
\rho({\cal F}_{\alpha ,k,l,t})\prod_{i\in \Z} \rho(h_{i})^{\gamma (k,l)_i}
\rho({\cal E} _{\beta ,k,l,t}) \in End \Vf \e(69)
$$
The map $\rho$ is a homomorphism of $\UA$ in $End\;\Vf$. It is
continuous in the $h-$adic topology. Hence $\rho$ defines a
representation of $\UA$ in $\Vf$, which is a highest weight
irreducible representation (since it is a highest weight
irreps with respect to the subalgebra $\Ua$).

Let $v=\sum_{k\in \Z_+}h^k v_k,\;\; v_k \in \Vh$ be an arbitrary
element from $\Vf$. Since each $v_k$ is a finite linear combination
of $C$-vectors, for any $k\in \Z_+$ there exists 
an integer $N_k>1$ such
that $v_k\in \Vf_{N_k}$.  Then from (59) and 
(61) one concludes
that the following properties hold:
$$
\eqalignno{
& 1. \q \rho(U_h(n_+)_{\b}) v_k=0, \;\; if \;\; S(\b) 
     \not\subset (-r_{N_k},r_{N_k}), & (70a)\cr
& 2. \q \rho(U_h(n_-)_{-\a}) v_k=0, \;\; if \;\; S(\a) \subset (-\infty,-r_{N_k}]
     \;\; or \;\; S(\a) \subset [r_{N_k},\infty), & (70b) \cr
& 3. \q \rho(h_i)v_k=0, \;\; if \;\; |i|\ge r_{N_k}.& (70c)\cr
}
$$
Therefore each finite signature representation of $\UA$ in $\Vf$
is a restricted representation (see {\it Definition 4.1} in Ref. [8]).

Let us mention in conclusion that the $\Ua$ modules $\V$, which
are not of a 
finite $(m,n)$ signature and for which $\xi_0\neq M_m,\;\xi_1\neq M_n,$ 
cannot be turned into $\UA$ modules.
In order to see this consider the transformation of the highest
weight vector $(\hat M)$ under the action of the operator
$I=\sum_{i\in \Z} h_i \in \UA$. From (20) and (41) one obtains:
$$
\r(I)(\hat M)=\left(\sum_{i=-\infty}^0 (M_i-\xi_0)+
\sum_{i=1}^\infty (M_i-\xi_1)\right)(\hat M). 
\e(71)
$$
The r.h.s. of 
(71) is not divergent only in the
finite-signature modules $\Vf$.

We were dealing with algebras and modules over $\Ch$.  The
algebra $\UA$ is however well defined also for $h$ being a
fixed complex number $h_c$, such that $h_c\notin i\pi \Q$,
namely
in the case 
$q=e^{h_c}$ is not a root of 1 $[8].$ In that case the
representations we have obtained remain highest weight irreps of
$U_{h_c}(A_\infty)$ realized in infinite-dimensional complex
linear spaces.

\vskip 18pt
\noindent
{\it Acknowledgments.}
We are grateful to Prof. Randjbar-Daemi for the kind hospitality
at the High Energy Section of ICTP.  One of us (N.I.S.) is
grateful to Prof. M.D. Gould for the invitation to work in his
group at the Department of Mathematics in University of
Queensland. The work was supported by the Australian Research
Council and by the Grant $\Phi-416$ of the Bulgarian Foundation
for Scientific Research.

\bigskip\n
{\bf References}

\vskip 12pt
{\settabs \+  $^{11}\;\;$ & I. Patera, T. D. Palev, Theoretical
   interpretation of the experiments on the elastic \cr

\+ $\;\;1. \;\;$ & Date, E., Jimbo, M., Kashiwara M., Miwa, T.:
           Operator Approach to the Kadomtsev-Petviashvili \cr
\+       & Equation - Transformation Groups for Soliton Equations III.		   
	     J. Phys. Soc. Japan {\bf 50}, 3806-3818 (1981).\cr

\+ $\;\;2.\;\;$ & Date, E., Jimbo, M., Kashiwara, M., Miwa, T.: 
           Transformation groups for soliton equations - Euclidean   \cr
\+       & Lie algebras and reduction of the KP hierarchy. Publ. RIMS 
           Kyoto Univ. {\bf 18}, 1077-1110 (1982).\cr

\+ $\;\;3.\;\;$ & Drinfeld, V.:  {\it Quantum groups.} ICM proceedings,
            Berkeley 798-820 (1986).\cr

\+ $\;\;4.\;\;$ & Goddard, P., Olive,D.: Kac-Moody and Virasoro
            Algebras in Relation to Quantum Physics. Int. J. \cr
\+&      Mod. Phys. {\bf A 1}, 303-414 (1986).\cr

\+ $\;\;5.\;\;$ & Haldane, F.D.M.: "Fractional statistics" in arbitrary
    dimensions: a generalization of the Pauli principle. \cr
\+ & Phys. Rev. Lett. {\bf 67}, 937-940 (1991).\cr

\+ $\;\;6.\;\;$ & Kac, V.G., Peterson, D.H.: Spin and wedge
         representations of infinite-dimensional Lie algebras and \cr
\+      &   groups. 
           Proc. Natl. Acad. Sci. USA {\bf 78}, 3308-3312 (1981). \cr

\+ $\;\;7.\;\;$ & Karabali, D., Nair,V.P.: Many-body states and
        operator algebra for exclusion statistics. Nucl. Phys. \cr
\+        &  {\bf B 438}, 551-560 (1995). \cr

\+ $\;\;8.\;\;$ & Levendorskii, S., Soibelman, Y.: Quantum group $A_\infty .$ 
             Commun. Math. Phys.
            {\bf 140}, 399-414 (1991).\cr

\+ $\;\;9.\;\;$ & Palev, T.D.: Lie algebraical aspects of the quantum
       statistics. Thesis, Institute of Nuclear Research \cr 
\+ &   and Nuclear Energy (1976),  Sofia; 
   Palev, T.D.: Lie algebraic aspects of quantum statistics. Unitary \cr

\+    & quantization (A-quantization). Preprint JINR E17-10550 
     (1977) and hep-th/9705032. \cr

\+ $10.\;\;$ & Palev, T.D.: Lie superalgebras, infinite-dimensional algebras
and quantum statistics. Rep. Math. \cr 
\+ & Phys. {\bf 31}, 241-262 (1992).\cr

\+ $11.\;\;$ & Palev, T.D.: Representations with a highest weight of the Lie
   algebra $a_\infty $. Funkt. Anal. Prilozh. \cr
\+ & {\bf 24}, No 3, 88-89 (1990) (in Russian); Funct. Anal. Appl. 
{\bf 24}, 250-251 (1990) (English translation).\cr

\+ $12.\;\;$ & Palev, T.D.: Highest weight irreducible unitarizable
   representations of the Lie algebras of infinite \cr
\+ & matrices. I. The algebra
$A_\infty $. Journ. Math. Phys. {\bf 31}, 1078-1084 (1990).\cr

\+ $13.\;\;$ & Palev, T.D., Stoilova, N.I.: Many-body Wigner quantum
    systems.  Journ. Math. Phys. {\bf 38},  2506-2523 \cr
\+ &(1997) and hep-th/9606011. \cr

\+ $14.\;\;$ & Palev, T.D., Stoilova, N.I.: Highest weight representations
of the quantum algebra $U_h(gl_\infty ).$  Preprint \cr 
\+   & ICTP IC/97/29
             (1997) and q-alg/9704001. \cr

\end